\documentclass{amsart}

\usepackage{graphicx}

\title{Cycles and Commutative Algebra}
\author{Paul C. Roberts}
\date{\today}

\newtheorem{thm}{Theorem}
\newtheorem{prop}{Proposition}

\newtheorem{defn}{Definition}

\newcommand{\p}{\mathfrak p}
\newcommand{\q}{\mathfrak q}

\newcommand{\m}{\mathfrak m}

\newcommand{\Px}{{\mathbb P}}
\newcommand{\cch}{\mbox{ch}}
\newcommand{\rank}{\mbox{rank}}
\newcommand{\length}{\mbox{length}}
\newcommand{\Tor}{\mbox{Tor}}
\newcommand{\Ext}{\mbox{Ext}}
\newcommand{\dime}{\mbox{dim}}

\newcommand{\depth}{\mbox{depth}}
\newcommand{\Spec}{\mbox{Spec}}
\newcommand{\Proj}{\mbox{Proj}}
\newcommand{\divv}{\mbox{div}}
\newcommand{\degree}{\mbox{degree}}

\begin{document}
\maketitle

The aim of this article is to describe several applications of the theory
of cycles 
to problems in Commutative Algebra.  The main topic is
the use of the theory of
 local Chern characters defined in the Chow group of a ring to answer some
questions on modules of finite homological dimension and to clarify others.

In the first section, we describe the origins of these problems in
Intersection Theory, and in particular Serre's definition of intersection
multiplicities and the conjectures which arose from this definition. We
next give an overview of the main facts on projective
dimension and related areas of Commutative Algebra. 
We then present an outline of the main properties of the Chow group and
 local Chern characters, together with 
their relation to modules of finite projective dimension.  Finally,
in the last two sections we
describe applications of the  theory and present the current
state of research in this area.

This paper is based on talks given at the conference on 
cycles in Morelia, Mexico in June 2003.

\section{Background: Serre's conjectures on intersection multiplicities}

Most of the questions we discuss in this paper originated from Serre's algebraic definition of
intersection multiplicities.  There have been many definitions of intersection
multiplicity proposed in Geometry and Topology, each with its own range of application
and set of assumptions.  In any of them the basic idea is to define
the order of tangency of two subspaces meeting at a point in such a 
way that certain natural conditions hold.  Serre's definition is
purely algebraic and will satisfy conditions listed below.

We start with a very simple example, that of the intersection of two curves in the
plane, and, in fact a simple case  where one of the curves is the $x$-axis.  
Let the other curve be defined by the equation $y=x^3-x^2$. 

It is clear (from Calculus, for example) that there are two points
of intersection, the origin $(0,0)$ and the point  $(1,0)$, the first with multiplicity
2 and the second with multiplicity 1.  The basic idea behind early algebraic definitions
of intersection multiplicity   is that it should be determined by the  dimension of the
vector space obtained by dividing the polynomial ring
$k[x,y]$ by the ideal generated by the polynomials defining the curves.  In this
case the polynomials are $y$ and $y-x^3+x^2$, and the quotient $k[x,y]/(y,y-x^3+x^2)$
is isomorphic to   
$k[x]/(x^3-x^2)$, which (by the Euclidean division algorithm) has dimension 3.  This
number gives  the total number of intersections counted with the appropriate
multiplicities.  To obtain the multiplicity at each point, we replace the
ring
$k[x,y]$ by its localization at that point.  For example, at the point $(0,0)$ in this
example, if $A$ denotes the localization of the polynomial ring at the ideal $(x,y)$, 
 it is a simple exercise to show that 
the dimension of $A/(y,y-x^3+x^2)$ is two.

This example suggests that we can define mutiplicities in general as follows:
we take the local ring at a point, which we denote
$A$, take the ideals defining the
two subvarieties near the point, say $I$ and $J$, and define the intersection
multiplicity as the vector space dimension $\dime_k(A/(I+J)).$  We remark that
the fact this point is an isolated point of the intersection assures that
this dimension is finite.

The problem with this definition is that it lacks some of the properties
required by intersection multiplicities; in particular,
it does not satisy B\'ezout's Theorem.  B\'ezout's theorem states that if we have two
subvarieties of projective space, and if they are of complementary dimension and
intersect in a finite number of points, then the number of points of intersection counted with
multiplicities is the product of the degrees of the subvarieties.
(For example, even though the above example is not in
projective space, the theorem does hold and both the sum of multiplicities and the
product of the degrees of the defining equations  are equal to 3.)

We now give Serre's definition of intersection multiplicities, which does not have
this drawback.  To improve flexibility, the definition is given in terms of modules
$M$ and $N$; the case of subvarieties is the case in which $M=A/I$ and $N=A/J$ as
above.

\begin{defn}
Let $A$ be a regular local ring, and let $M$ and $N$ be finitely generated modules such
that $M\otimes_A N$ is a module of finite length.  Then the 
intersection mutiplicity of $M$ and $N$ is 
$$\chi(M,N)=\sum_i (-1)^i\mathop{\rm {length}}( \mathop{\rm {Tor}}\nolimits_i^A(M,N)).$$
\end{defn}

It may seem surprising that a rather complicated definition in terms of Euler
characteristics such as this would satisfy a simple formula like B\'ezout's Theorem.  
We digress a little and explain why it does, and for this we recall a little more
history. The introduction of homological methods in Commutative Algebra  goes
back (at least) to Hilbert.  He introduced graded resolutions of modules over graded
rings,  and he  also introduced  what are now called Hilbert polynomials.  Let $A$ be
a graded ring such that $A_0$ is a field and $A$ is generated over $A_0$ by $A_1$ (for
example, a graded polynomial ring).  If
$M$ is a finitely generated graded module over
$A$, then we let $H_M(n)$  equal  the length (that is, the vector space dimension over the
field $A_0$) of the
graded piece
$M_n$ for all integers
$n$. Hilbert showed that
there is a polynomial $P_M(n)$ such that $P_M(n)=H_M(n)$ for
large $n$; this polynomial is now called the Hilbert polynomial of $M$.    Suppose that
$$0\to F_k\to \cdots\to F_1\to F_0\to M\to 0$$ is a free resolution of $M$ by graded
free modules; that is, each $F_i$ is a finite direct sum of modules of the form
$A[n_{ij}]$, where for any integer $k$, $A[k]$ denotes the graded module $A$ with
degrees shifted by $k$ so that $A[k]_m=A_{m+k}$ for all $m$.  Hilbert showed that over
a polynomial ring every module has such a resolution (The Hilbert Syzygy Theorem). 
Furthermore, since Hilbert polynomials are additive on short exact sequences, this gives a way to
compute the Hilbert polynomial:
$$P_M(n)=\sum_{i=0}^k(-1)^iP_{F_i}(n)=\sum_{i=0}^k\sum_{j}(-1)^iP_{A}(n+n_{ij}).$$
Finally, if the graded ideal $I$ defines a closed subscheme $\Proj(A/I)$ of
projective space, the degree of this subscheme is simply the
leading coefficient of the Hilbert polynomial $P_{A/I}(n)$ multiplied by $d!$,
where $d$ is the dimension of $\Proj(A/I)$. 
The point of this discussion is that both intersection multiplicities and degrees of
subvarieties
are defined by alternating sums of lengths of modules derived from free resolutions; combining all
this with the fact that everything commutes with localization  is the idea behind B\'ezout's Theorem,
which we leave as an (admittedly nontrivial) exercise.

We now return to the discussion of intersection multiplicities.  
The definition requires two conditions: that $\Tor_i(M,N)$ have
finite length for each $i$ and that it be zero for large $i$.  The
first condition follows from the assumption that $M\otimes N$ has finite length.  The
second condition is a modern version of the Hilbert Syzygy Theorem, that for a
regular local ring, the projective dimension of any module is finite.

The intersection multiplicity has the following properties:

\begin{enumerate}
\item It satisfies B\'ezout's Theorem.
\item If $M=A/I$ and $N=A/J$, where $I$ and $J$ are ideals defined by
smooth subschemes intersecting transversally, then $\chi(M,N)=1$.
\end{enumerate}

Serre also made the following conjectures, one of which is still open.

Let $A$ be a regular local ring, and let $M$ and $N$ be finitely generated $A$-modules such that $M\otimes_AN$ has finite length.

\begin{enumerate}
\item $\dime(M)+\dime(N)\le \dime(A).$
\item $\chi(M,N)\ge 0.$
\item $\chi(M,N)>0$ if and only if $\dime(M)+\dime(N)=\dime(A).$
\end{enumerate}

Serre proved the first in general and proved the second and third  by reduction to the
diagonal in the equicharacteristic case.  The second and third are equivalent to
\begin{enumerate}
\item [$2'.$] If $\dime(M)+\dime(N)<\dime(A)$, then $\chi(M,N)=0.$
\item [$3'.$] If $\dime(M)+\dime(N)=\dime(A)$, the $\chi(M,N)>0$.
\end{enumerate}

These conjectures and generalizations of them have been among the
central problems in Commutative Algebra since that time.  The first
to be proven was $2'$, by myself (\cite{r-vanish},\cite{r-it}) and Gillet and Soul\'e
(\cite{gil-sou-cr},\cite{gil-sou}) around 1985.  I will discuss my proof of this briefly later, when
I talk about local Chern characters.  Conjecture 2
was proven by Gabber around 1996; this was the topic of the talks I
gave at Guanajuato in 1997 \cite{R?}.  

The other questions that we discuss concern to the extent to which these 
properties can be proven
without the assumption that $A$ is regular, and this may be summed
up by asking how many of the results follow from the assumption that
the modules involved have
finite projective dimension in the nonregular case.  This will be discussed in
 the last sections of this paper.  As we shall see, all of these topics use 
the theory of  cycles in
fundamental ways.

\section{Assorted facts about modules of finite projective dimension}

In this section we recall some basic facts from Commutative Algebra that
relate to modules of finite projective dimension. Throughout this
section $A$ will denote a local ring with maximal ideal $\m$, and modules will be
assumed to be finitely generated.  Most ofhe results here can
be found in Matsumura \cite{matsu}, Chapters 6 and 7.

Any module $M$ has a free resolution:
$$\cdots \to F_k\to F_{k-1}\to\cdots\to F_1\to F_0\to M\to 0.$$
In general, the resolution goes on forever, and $M$ is said to have
{\em finite projective dimension} if such a resolution can be chosen
so that there exists an $n$ with $F_i=0$ for $i>n$.  The smallest
$n$ for which this is possible is the {\em projective dimension}
of $M$, denoted pd$(M)$ (we are ignoring the distinction between projective and free modules
since we assume the ring $A$ is local, so that all projective modules are free).

 If $M$ is finitely
generated, it has a minimal set of generators $m_1,\ldots m_s$ that form
a basis of the vector space $M/\m M$ over the field $A/\m$.  These have
the property that every relation $a_1m_1+\cdots +a_sm_s=0$ has all
coefficients $a_i$ in $\m$, and any two minimal sets of generators have the same number of elements
and  can be written in terms of each other  by an invertible matrix.
A {\it minimal free resolution} of $M$ is obtained
by taking  a minimal set of
generators at each step; the result is a complex such that the matrix defining
each map of free modules has  coefficients in
$\m$. Any two minimal free resoutions of $M$ are isomorphic as complexes. Furthermore, the last
nonzero
$n$ such that $F_n\not=0$ in a minimal free resolution is the projective dimension of $M$.

We next recall some facts about depth and its relation to modules
of finite projective dimension.  If $M$ is a nonzero finitely generated
module, a sequence of elements of $\m$, $x_1,\ldots, x_s$ is {\em $M$-regular}
if $x_i$ is not a zero divisor modulo $(x_1,\ldots, x_{i-1})M$ for
$i=1,\ldots, s$.  A classical example of an $M$-regular
sequence is where $M=A=$ a power series ring on the variables $x_1,\ldots, x_s$.
The {\em depth} of $M$ is the length of the longest regular sequence;
it can be shown that every maximal $M$-regular sequence has
the same length and that this is at most equal to the 
dimension of the module.
If the depth is equal to the dimension, the module is {\em Cohen-Macaulay};
a ring is Cohen-Macaulay if it is Cohen-Macaulay as a module.

A local ring is {\em regular} if its maximal ideal is generated by $d$ 
elements, 
where $d$ is the dimension of 
$A$.  If so, and if $x_1,\ldots, x_d$ are a minimal
set of generators for the maximal ideal, then $x_1,\ldots, x_d$ form
a regular sequence and $A$ is Cohen-Macaulay.
 
A module $M$ has depth zero if and only if there is a nonzero element of $M$ that is
annihilated by
the maximal ideal $\m$.

We have the following facts:

\begin{enumerate}
\item The Auslander-Buchsbaum-Serre Theorem: the local ring $A$ is
regular if and only if every $A$-module has finite projective dimension.
Part of this theorem generalizes the Hilbert Syzygy Theorem as mentioned above.

\item The Auslander-Buchsbaum formula: if $M$ has finite projective dimension then
$$\mbox{pd}(M)+\depth(M)=\depth(A).$$
We remark that if $A$ has depth zero this theorem implies that the only modules
of finite projective dimension are free, which is in fact
quite easy to show: if $M$ has finite projective dimension and is not free, its
minimal free resolution ends with
$$0\to F_n\to F_{n-1}\to $$
for some $n>0$.  Since the resolution is minimal the matrix defining
this map has entries in $\m$, but since $A$ has depth zero there
is an element of $F_n$ annihilated by $\m$, so this map cannot be injective,
giving a contradiction.

\item The Acyclicity Lemma: Let 
$$0\to M_n\to M_{n-1}\to \cdots \to M_k\to \cdots
$$
be a complex of  modules with homology of finite length.  Assume that each module
$M_i$ has depth at least $r$. Then
$H_i(F_\bullet)=0$ for $n \ge i \ge n-r+1.$

\end{enumerate}

We note a consequence of these formulas for a special case of a generalized form of Serre's
positivity conjecture over a Cohen-Macaulay ring $A$.  Suppose that $M$ and $N$ are Cohen-Macaulay
modules such that $M$ has finite projective dimension, $M\otimes_A N$ has
finite length, and $\dime(M)+\dime(N)=\dime(A)$.  Then
if $$0\to F_n\to F_{n-1}\to\cdots\to F_0\to M\to 0$$
is a minimal free resolution, we have $n=\dime(A)-\dime(M)$ by the
Auslander-Buchsbaum formula, since both $A$ and $M$ are Cohen-Macaulay.
If we tensor with $N$, we get a complex of length $n$. Since we are assuming
that $\dime(M)+\dime(N)=\dime(A)$, $n$ is the 
 dimension of $N$, which is the depth since $N$ is also assumed to be Cohen-Macaulay. 
Thus higher Tors are zero by the Acyclicity Lemma, and
$$\chi(M,N)=\length(M\otimes_A N)>0.$$

It follows from this argument that  if finitely generated Cohen-Macaulay modules
exist with given support, Serre's positivity conjecture  (and several other conjectures)
could be proven. However, the existence of such modules, called ``small Cohen-Macaulay
modules", is known in very few cases.

In the final section we will return to these questions on multiplicity
properties of finite projective dimension, after the introduction
of the Chow group and a discussion of local Chern characters.

\section{The Chow group and local Chern characters}

The use of the Chow group and local Chern characters has made it
possible to prove some of the conjectures described above, and
it also led to a better understanding and strengthening of
some of the results.

For any scheme $X$ of finite type over a regular scheme, the Chow
group
$CH_*(X)$ is defined to by the group of cycles modulo rational equivalence.
There are two cases of special interest here, and we describe these
in more detail.

If $A$ is a Noetherian ring, then for each integer $i\ge 0$, we let
$Z_i(A)$ be the free abelian group on the prime ideals $\p$ of
$A$ such that the dimension of $A/\p$ is equal to $i$.  For each
prime ideal $\q$ with $\dime(A/\q)=i+1$ and for each $f\not=0$ in
$A/\q$, define $$\divv(f,A/\q)=\sum\length(A/(\q,f))_{\p}[A/\p],$$
where the sum is taken over all prime ideals $\p$ such that 
dim$(A/\p)=i$.  
The component of dimension $i$ of the Chow group, $CH_i(A)$, is then the quotient of
$Z_i(A)$ by the subgroup generated by all $\divv(f,A/\q)$ for all such $\q$
and $f$.

If $A$ has dimension $d$, then $CH_d(A)$ is the free abelian group on
the components of $\Spec(A)$ of dimension $d$.  If $A$ is an
integrally closed domain of dimension $d$, then $CH_{d-1}(A)$ is
the ideal class group of $A$.

The other case is where $A$ is a graded module over a field and $X$ is
the associated projective scheme.  In this case the description is
similar except for two major differences:
\begin{enumerate}
\item  $Z_i(X)$ is generated by graded prime
ideals with
$\dime(A/\p)=i+1$ (so that the projective subscheme defined by $A/\p$ has
dimension $i$).

\item  The relation of rational equivalence is defined by setting  $\divv(f,A/\q)=0$,
where $\q$ is a graded prime ideal and $f$ is a
quotient of two homogeneous polynomials of the same degree. Thus $f=g/h$, and
$\divv(f,A/\q)=\divv(g,A/\q)-\divv(h,A/\q)$ is zero in $CH_*(X)$ (but neither
$\divv(g,A/\q)$ nor $\divv(h,A/\q)$ is necessarily zero in $CH_*(X))$.

\end{enumerate}

 Note that in the case in which $A$ is graded and $X=\Proj(A)$
there is a map from $CH_i(X)$ to $CH_{i+1}(A)$ induced by the inclusion
of the set of graded prime ideals into the set af all prime ideals of $A$.  We will
state the precise relation between these two Chow groups in a later section.

For $X=\Proj(A)$ there is an important operator called the {\it hyperplane section}
on $CH_*(X)$; we denote this operator $h$.  It is defined as the
map from $CH_i(X)$ to $CH_{i-1}(X)$ that sends a generator  $[A/\p]$
to $\divv(A/\p,x)$, where $x$ is any homogeneous element of $A$ of
degree 1 that is not in $\p$.  It is easy to check that this
definition is well defined up to rational equivalence.

We now discuss how to  reinterpret intersection multiplicities in
terms of actions on the Chow group.
  The main tool that connects
the two theories is the use of local Chern characters.

We recall that we had defined intersection multiplicities of
two modules $M$ and $N$ over a regular local ring in terms of Tor modules. 
This means that we compute $\chi(M,N)$ by taking a free
resolution of one, say $F_\bullet\to M\to 0$, tensoring with the other module $N$, and taking the
alternating sum of lengths of homology modules.  The
assumptions necessary so that this multiplicity is defined
are that there are only a finite number of nonzero Tor modules and the all Tor modules have finite
length.  We generally assume in what follows that $M$ has finite projective dimension and that
$M\otimes_A N$ has finite length, so that both these conditions hold.

   We first generalize this definition  to perfect complexes, where we say that
a complex is {\em perfect} if it is of the form
$$0\to F_k\to
\cdots
\to F_0\to 0$$ where each $F_i$ is free.  A special case is a minimal free resolution of a module
of finite projective dimension. 

 If $F_\bullet$ is a perfect complex and  $N$ is a module such that $F_\bullet \otimes N$ 
has homology of finite length, we define 
$$\chi(F_\bullet,N)=\sum(-1)^i\length(H_i(F_\bullet\otimes N)).$$

We now outline the theory of local Chern
characters on an scheme $X$.  Let $F_\bullet$ be a perfect complex with support $Z$.  Then
the local Chern character $\cch(F_\bullet)$ is a map, for
each subscheme $Y$ of $X$, from $CH_*(Y)$ to $CH_*(Y\cap Z)\otimes_{\mathbb Z}{\mathbb Q}$,
which we denote $CH_*(Y\cap Z)_{\mathbb Q}$.
More precisely, for each $i$ there is a component
$\cch_i$ from $CH_j(Y)$ to $CH_{j-i}(Y\cap Z)_{\mathbb Q}$ for each $j$. These maps satisfy
a number of functoriality conditions, of which we list the
ones we need below.  We will not give the definition of local Chern characters,
which is very involved, but we will give some examples to give an
idea of what they are.  For a complete description we refer to Fulton \cite{fulton}
and Roberts \cite{ccm}.

{\bf Examples:} Let $F_\bullet$ be a perfect complex with support $Z$.  We give
some examples of the local Chern character $\cch_i(F_\bullet)$.
\begin{enumerate}
\item The case $i=0$:  this is multiplication by the alternating sum of the ranks of
$F_i$.  That is, if $r=\rank(F_0)-\rank(F_1)+\cdots +(-1)^k\rank(F_k)$, then
$\cch_0(F_\bullet)_{\mathbb Q}$ is the map from $CH_j(Y)$ to $CH_j(Y\cap Z)_{\mathbb Q}$ defined by
multiplication by $r$.  Note that if $Z$ is a proper subset of $\Spec(A)$,
we must have $r=0$, so we get the zero map.  If $r$ is nonzero we have $Y\cap Z=Y$ and this map
is multiplication by $r$ from $CH_*(Y)$ to $CH_*(Y)_{\mathbb Q}$.

\item   Suppose $F_\bullet$ is $0\to F_1\stackrel{x}{\to}F_0\to 0,$
where $F_0=F_1=A$, $x$ is an element of $A$, and $A$ is an integral domain.  Then
we have the intersection with the divisor defined by $x$, which takes a generator
$[A/\p]$ of the Chow group of $A$ to 
class we denoted $\divv(x,A/\p)$ in the Chow group of the  quotient $A/xA$.  Note
that although $\divv(x,A/\p)$ is zero in the Chow group of $A$ by definition,
it is not necessarily zero in the Chow group of $A/xA$.    In this case 
$\cch_0(F_\bullet)=0,$ $\cch_1(F_\bullet)$ is intersection with the divisor
defined by $x$, and $\cch_i(F_\bullet)$ is zero for $i\ge 2$.

\item 
More generally, if $F_\bullet$ is of the form
$$0\to F_1\to F_0\to 0,$$
where $F_1$ and $F_0$ have  equal rank, the local Chern
character  is intersection with
the determinant of the matrix defining the map from $F_1$ to $F_0$.  This example was generalized
to an intersection operator of codimension 1 for longer complexes  as a generalization of the 
 MacRae invariant 
by Foxby
\cite{fox} and was one of the  precursors of the general theory.

\item
We stated that if we have a complex of the form $0\to A\stackrel{x}{\to}
A\to 0$, where $A$ is a Noetherian ring, then the local Chern character is
intersection with the divisor defined by $x$.  For projective schemes,
the situation is different.  Let $A$ be a standard graded
ring and let $X=\Proj(A)$.  We consider the complex
$$0\to A[-n]\stackrel{x}{\to}A\to 0$$
where $x$ is an element of degree $n$.  In this case the local Chern character
is multiplication by the element $1-e^{-nh}$ where $h$ is the hyperplane
section.  We note that this comes out to
$$1-(1-nh+{{n^2}\over {2!}}h^2-\cdots)=nh-{{n^2}\over {2!}}h^2+{{n^3}\over {3!}}h^3\cdots.$$
More precisely, this expression defines the image of $\cch(F_\bullet)([X])$ in
the Chow group of $X$ over $\mathbb Q$; unlike the affine case, this image is not zero.  The
element $\cch(F_\bullet)([X])$ itself in the Chow group of $\Proj(A/xA)$ is given by
$$(1-{{n}\over {2!}}h+{{n^2}\over {3!}}h^2-\cdots)([\Proj(A/xA)]).$$
The reason that this reduces to  intersection with $x$ in the affine case is that in that case the
divisor intersected with itself is zero in the Chow group.

\end{enumerate}

We now describe the basic properties of local Chern characters that are used in applications to
Commutative Algebra. 

In many of the cases we consider, the perfect complex $F_\bullet$ has homology of
finite length, so that its support is $\Spec(A/\m)$, and the local Chern character
defines a map from $CH_*(A)$ to $CH_*(A/\m)_{\mathbb Q}=CH_0(A/\m)_{\mathbb Q}
\cong {\mathbb Q}$. 
Thus for any element
$\eta$ of the Chow group, there is a rational number $\cch(F_\bullet)(\eta)$.  On the other
hand, for any module $M$ we have an integer $\chi(F_\bullet\otimes M)$.  One
of the main theorems in this subject, the local Riemann-Roch theorem, states that
there is a map from the Grothendieck group of the category of bounded
complexes of finitely generated modules to
the Chow group that relates these invariants.  More precisely, if $M_\bullet $ is a 
bounded complex of finitely
generated modules, there is a class $\tau(M_\bullet)$ in the Chow group of the support of $M_\bullet$
such that the following theorem holds. 

\begin{thm} (The Local Riemann-Roch formula).  Let $F_\bullet$ be a perfect complex and let
$M_\bullet$ be a bounded complex of finitely generated modules.  Then
$$\tau(F_\bullet\otimes_A M_\bullet)=\cch(F_\bullet)(\tau(M_\bullet)).$$
\end{thm}

Furthermore, if $M_\bullet$ has homology of finite length, we have 
$$\tau(M_\bullet)=\chi(M_\bullet).$$
under the identification of $CH_0(A/\m)_{\mathbb Q}$ with $\mathbb Q$.  A particular case is when
$M=A$; there is a class $\tau(A)$ such that
$$\chi(F_\bullet)=\chi(F_\bullet\otimes_A A)=\cch(F_\bullet)(\tau(A)).$$

The dimension $d$ component of $\tau(A)$ is the class of the components of dimension $d$; that is,
it is $\sum n_{\p}[A/\p]$, where the sum is over primes $\p$ with $\dime(A/\p)=d$ and $n_{\p}$ is
the length of the localization $A_{\p}$. 

The local Chern character also satisfies the following properties:

\begin{enumerate}
\item (Additivity) If $0\to F'_\bullet\to F_\bullet\to F''_\bullet\to 0$
is a short exact sequence of perfect complexes, then
$$\cch(F_\bullet) = \cch(F'_\bullet) + \cch(F''_\bullet).$$
\item (Multiplicativity) If $F_\bullet$ and $G_\bullet$ are perfect complexes,
then $$\cch(F_\bullet\otimes G_\bullet)=\cch(F_\bullet)\cch(G_\bullet)$$
in the Chow group of the intersection of the supports of $F_\bullet$ and $G_\bullet$.
\item (Commutativity) If $F_\bullet$ and $G_\bullet$ are perfect complexes, then
for any $i$ and $j$ we have
$$\cch_i(F_\bullet)\cch_j(G_\bullet)=\cch_j(G_\bullet)\cch_i(F_\bullet).$$
\item (Additivity of $\tau$) If $0\to M'_\bullet\to M_\bullet\to M''_\bullet\to 0$
is a short exact sequence of bounded complexes of finitely generated modules, then
$$\tau(M)=\tau(M')+\tau(M'').$$
\end{enumerate}

\section{An extended example: graded modules}

In 1975, Peskine and Szpiro \cite{ps-syz-mult} proved that 
virtually all of the multiplicity properties we 
have been discussing
hold for graded modules over graded rings.  This theorem relied on a formula for Hilbert polynomials
that follows from the Taylor expansion.  This result was one of the motivations for the development
 of 
the method of local Chern characters and, when stated in the right way, can be thought of
as a special case of this theory.  We present this example here.

Let $A$ be a graded ring such that $A_0$ is a field (in fact, $A_0$ can be taken to be 
any
Artinian ring) and $A$ is generated over $A_0$ by $A_1$.  Let $X=\Proj(A)$,
and, as before, let $h$ be the hyperplane section on the Chow group of $X$. 

Let $M$ be a graded module of finite projective dimension.  Then $M$ has a finite resolution
by graded free modules 
$$0\to F_k\to \cdots\to F_1\to F_0\to 0,$$
where we have $F_i=\oplus_j A[n_{ij}]$ for some integers $n_{ij}$.  Since local Chern characters are
additive, it is easy to compute the image of the action of the local Chern character of
$F_\bullet$.  In fact, it is simply the alternating sum of the Chern characters
defined by the $A[n_{ij}]$, and the Chern character of $A[n_{ij}]$ is $e^{n_{ij}h}$, as we showed in
the previous section.  Combining these,and expanding $e^{n_{ij}h}$ in a power series, we have
$$\cch(F_\bullet)=\sum(-1)^i\cch(F_i)=\sum_{i,j}(-1)^ie^{n_{ij}h}=\sum_{i,j}(-1)^i\left(\sum_k
{{n_{ij}^kh^k}\over{k!}}\right).$$ Hence we can write
$$\cch(F_\bullet)=\sum \rho_k(F_\bullet)h^k,$$
where $$\rho_k(F_\bullet)=\sum_{i,j}(-1)^i{{{n_{ij}^k}\over{k!}}}.$$  

Thus the local Chern character of the complex $F_\bullet$   can be
explicitly described by the above formula.  The definition of the $\tau(M_\bullet)$  is a little less
straightforward.  Assume that $X$ is embedded in $\Px^N$.  We can
map the Chow group of $X$ to the Chow group of $\Px^N$.  The Chow group of 
$\Px^N$ over $\mathbb Q$ can be identified with the ring ${\mathbb Q}[h]/h^{N+1}$,
where the element $h^k$ is identified with $h^k([\Px^N])$.  The element
we define is in the ring ${\mathbb Q}[t][h]/h^{N+1}$, where $t$ is an indeterminate.

For this application we use the following version of the Hilbert polynomial for graded modules
and complexes.
If $M$ is a graded module, we let 
$$P_M(t)=\sum_{k\le t}\dime_{A_0}M_k,$$
where $M_k$ is the graded piece of degree $k$.
We note that with this definition the degree of the Hilbert polynomial of $M$ is equal to the
Krull dimension of $M$, and the Hilbert polynomial is zero if and only if $M=0$.  If
$M_\bullet$ is a bounded complex of finitely generated modules, we let
$$P_{M_\bullet}(t)=\sum_i(-1)^iP_{M_i}(t)=\sum_i(-1)^iP_{H_i(M_\bullet)}(t);$$
here $M_i$ denotes the $i$th module in the complex $M_\bullet$.

  Let
$M_\bullet$ be a bounded complex of  finitely generated graded $A$-modules.  We define
$$\tau({M_\bullet})=P_{M_\bullet}(t)h^N+P'_{M_\bullet}(t)h^{N-1}+\cdots
+P_{M_\bullet}^{(k)}(t)h^{N-k}+\cdots$$ where $P_{M_\bullet}^{(k)}(t)$ denotes the $k$th 
derivative of the Hilbert polynomial $P_{M_\bullet}(t)$ of
 $M_\bullet$.

\begin{thm} With the above notation, we have
$$\tau(F_\bullet\otimes M_\bullet) = \cch(F_\bullet)\tau(M_\bullet).$$
\end{thm}
{\bf Proof.}  Since both sides of this equation are
additive on short exact sequences in both variables, we can
reduce the question to the case in which  $M_\bullet$ is one module $M$ and $F_\bullet$
is one module of the form 
$A[m]$ for
some integer $m$. Thus we must show that for a graded module $M$ we have
$$\tau(A[m]\otimes M)= \tau(M[n]) = \cch(A[m])\tau(M).$$
Both sides are poynomials in $h$ with coefficients in ${\mathbb Q}[t]$, and
the required equality for the coefficient of $h^{N-k}$ is

$$P_M^{(k)}(t+m)=\sum P_N^{(k+i)}(t)m^i/i!.$$ 
This is just the Taylor expansion of the $k$th derivative of $P_M(t)$.

\bigskip

The reason this formula implies the intersection conjectures for graded modules  is that the
coefficients
$\rho_k(F_\bullet)$ are simply rational numbers and the derivatives of the Hilbert polynomial
$P_N(t)$ are linearly independent polynomials.  We can thus deduce the following proposition.

\begin{prop}\label{rho}
Let notation be as above, and let 
$$q(F_\bullet)=\inf\{k|\rho_k(F_\bullet)\not=0\}.$$
Then
\begin{enumerate}
\item We have $\rho_k(F_\bullet)=0$ for all $k\le \degree(P_{M_\bullet}(t))$ if and
only if $P_{F_\bullet\otimes M_\bullet}(t)$ is the zero polynomial.
\item If there is a $k$ with $k\le \degree(P_{M_\bullet}(t))$ such that $\rho_k(F_\bullet)\not=0$,
then
$$q(F_\bullet)=\degree(P_{M_\bullet}(t))-\degree(P_{F_\bullet\otimes M_\bullet}(t)).$$
\item If not all $\rho_k$ with 
$k\le \degree(P_{M_\bullet}(t))$ are zero, 
the leading coefficients of $P_{M_\bullet}(t)$
and
$P_{F_\bullet\otimes M_\bullet}(t)$ have the same sign if and only if 
$\rho_{q(F_\bullet)}>0$. 
\end{enumerate}
\end{prop}

  We now apply this to the
situation in which
$F_\bullet$ is a free resolution of a module $M$ of finite projective dimension and $N$ is an
arbitrary finitely generated graded module.  

\begin{thm} Let $M$ be a nonzero graded module of finite projective dimension, and let $N$
be a graded module such that $M\otimes_A N$ has finite length.  Then 
\begin{enumerate}
\item $\dime(M)+\dime(N)\le \dime(A).$
\item $\chi(M,N)\ge 0.$
\item $\chi(M,N)>0$ if and only if $\dime(M)+\dime(N)=\dime(A).$
\end{enumerate}
\end{thm}  Let $F_\bullet$ be a graded free resolution of $M$. Applying Proposition \ref{rho}
to the case in which $M_\bullet=A$, we
obtain
$$q(F_\bullet)=\degree(P_A(t))-\degree(P_M(t))=\dime(A)-\dime(M).$$
Furthermore, since the leading coefficient of the Hilbert polynomial of a module is always
positive, we have $\rho_{q(F_\bullet)}>0.$

Applying Proposition \ref{rho} to the case $M_\bullet=N$, we have
$$q(F_\bullet)=\degree(P_N(t))-\degree(P_{F_\bullet\otimes N}(t)).$$
Since $M\otimes N$ has finite length, the polynomial $P_{F_\bullet\otimes N}(t)$ is constant
with value $\chi(M,N)$. Thus $ \degree(P_{F_\bullet\otimes N}(t))\le 0$, so
$$\dime(N) =\degree(P_N(t))\le q(F_\bullet)=\dime(A)-\dime(M).$$
Thus we have 
$$\dime(M)+\dime(N)\le \dime(N).$$
Furthermore, we have equality if and only if $q(F_\bullet)=\dime(N)$, and Proposition \ref{rho}
implies that this is true if and only if $\chi(M,N)$ is a nonzero constant.  Finally,
if we do have equality, since $\rho_{q(F_{\bullet})}>0$ it follows 
 that $\chi(M,N)>0$.

Thus the multiplicity conjectures hold for graded modules. 
They do not hold in this generality for arbitrary modules over a local ring.  While
there is an analogy between the methods of this section and the general theory
of local Chern characters, in general the local Chern character is an operator
on the Chow group and not simply a ratioanl number, and it is possiible for
$\cch_i(F_\bullet)$ to be nonzero for $i<\dime(A)-\dime(M).$
The examples discussed in the last section of the paper have this property.

\section{General vanishing theorems for multiplicities.}

In the last two sections we discuss vanishing theorems in various
situations and discuss several open problems.  In this section we consider the case
in which both $M$ and $N$ have finite projective dimension.  The
main theorem is the following.

\begin{thm}\label{vanishing} Suppose $M$ and $N$ are modules of finite projective dimension
over a ring of dimension $d$ such that $M\otimes N$ has finite length.  Suppose that
\begin{enumerate}
\item $\dime(M)+\dime(N)<d.$
\item $\tau(A)=[A]_d.$
\end{enumerate}
Then $\chi(M,N)=0$.
\end{thm}
{\bf Proof:} Since $M$ and $N$ have finite projective dimension, we
can find perfect complexes $F_\bullet$ $G_\bullet$ such that $F_\bullet\to M$ and
and $G_\bullet\to N$ are free resolutions.  Then $\Tor_i(M,N)$ is the $i$th 
homology of the complex $F_\bullet\otimes G_\bullet$, so we have
$$\chi(M,N)= \chi(F_\bullet\otimes G_\bullet).$$
Thus we have 
$$\chi(M,N)=\chi(F_\bullet\otimes G_\bullet)=
\cch(F_\bullet\otimes G_\bullet)(\tau(A))=
\cch(F_\bullet\otimes G_\bullet)([A_d])=$$
$$\cch(F_\bullet)\cch(G_\bullet)([A_d])=
\sum_{i+j=d}\cch_i(F_\bullet)\cch_j(G_\bullet)([A_d]),$$
where we have used the local Riemann-Roch formula, the hypothesis that 
$\tau(A)=[A]_d,$ and the multiplicativity of local Chern characters.
Suppose that
$j<\dime(A)-\dime(N)$, or, equivalently, that
$\dime(N)<d-j$.  Then there are no prime ideals
$\p$ in the support of $N$ with $\dime(A/\p)=d-j$, so $CH_{d-j}(\mbox{Supp}(N))=0$. Since
$\cch_j(G_\bullet)([A]_d)$ is an element of $CH_{d-j}(\mbox{Supp}(N))=0$, this
implies that $\cch_j(G_\bullet)([A]_d=0$, so that the term
$\cch_i(F_\bullet)\cch_j(G_\bullet)([A_d])$ in the above sum is zero.  Similarly, if
$i<\dime(A)-\dime(M)$, we have
$$\cch_i(F_\bullet)\cch_j(G_\bullet)([A_d])=\cch_j(G_\bullet)\cch_i(F_\bullet)([A_d])=0,$$
using the commutativity of local Chern characters.  
However, since $\dime(M)+\dime(N)<\dime(A)$, for
any
$i$ and
$j$ with
$i+j=d$, we must have $i<\dime(A)-\dime(M)$ or $j<\dime(A)-\dime(N)$.  Thus every term in the sum is
zero, so
$\chi(M,N)=0$.  
\bigskip

The class of rings for which the hypothesis that $\tau(A)=[A]_d$ includes complete intersections
and thus in particular regular local rings, so this result holds in the case
of Serre's original conjecture.  This
property has been studied in detail by Kurano \cite{kur-r}, and he has shown in particular
that it does not hold for general Gorenstein rings.

We remark also that  Jean Chan \cite{chan} has used this technique to show that
intersection multiplicities defined by $\Ext$ instead of $\Tor$ give the
same answer for complete intersections.

 One
of the main open questions in this field is whether Theorem \ref{vanishing}
holds without the assumption that $\tau(A)= [A]_d$.

\section{Modules of finite projective dimension with prescribed intersection properties.}

 In this section we assume that $M$ and $N$ are $A$-modules as above and that
$M$ has finite projective dimension but $N$ does not.  In this case, an
example of Dutta, Hochster, and McLaughlin \cite{dhm} showed that it is
possible for $\chi(M,N)$ to be nonzero (and in fact negative) if 
$\dime(M)+\dime(N)<\dime(A)$.  Recently, in joint work with Srinivas \cite{r-sri},
we showed that there are many such examples.  We describe the main
result here.

Assume that $R$ is a standard graded ring such that $R_0$ is a field, and assume that $X=\Proj(R)$ 
is smooth.  Assume also that the intersection pairing on the
Chow group $CH_*(X)$ is perfect (this is a very strong assumption, and
the main theorem of Roberts and Srinivas \cite{r-sri} is much more general, but this
case already includes many significant examples).  Since $X$ is a projective scheme,
the hyperplane section $h$ acts on the Chow group of $X$.  Let $A$ be the localization
of $R$ at the graded maximal ideal. It then follows from results of Kurano \cite{kurano}
that we have an isomorphism
$$CH_*(A)\cong CH_*(X)/hCH_*(X).$$
The main result of \cite{r-sri} in the case we are considering then states the following.

\begin{thm}\label{existence}  With notation as above, let $\eta$ be an element of $CH_i(X)$ such that
$h\cdot \eta=0$.  Then there is a perfect complex $F_\bullet$ over $A$ with homology of
finite length such that 
$$\chi(F_\bullet\otimes A/\p)=\eta\cdot [\Proj(A/\p)]$$
for all graded prime ideals $\p$, where $\eta\cdot [\Proj(A/\p)]$ denotes the
intersection pairing on $X$.  

If $A$ is Cohen-Macaulay, there is
 a module of finite length and finite projective dimension $M$ such that
$$\chi(M,A/\p)=\eta\cdot [\Proj(A/\p)]$$
for all graded prime ideals $\p$ with $\dime(A/\p)<\dime(A)$.
\end{thm}
\bigskip
We remark that if $A$ is not Cohen-Macaulay, it follows from the Intersection Theorem (see Peskine
and Szpiro \cite{ps-dpf} and Roberts \cite{r-it}) that
there are no modules of finite length and finite projective dimension, so the
first version of the theorem is the best possible.  

 This theorem provides nice examples when the Chow group of $X$ is easy to compute.  
We show how to prove the existence of an example in the situation used
in the paper of Dutta, Hochster, and McLaughlin mentioned above.  In this
case $R=k[X,Y,Z,W]/(XY-ZW)$, so $X=\Proj(R)=\Px^1\times \Px^1$.  Let $x$ be
a point of $X$, and let $Y$ and $Z$ denote
the subschemes $\{x\}\times \Px^1$ and $\Px^1\times \{x\}$ respectively.
Note that $h$ is intersection with the cycle $[Y]+[Z]$ and that $[Y]\cdot [Y]$ and
$[Z]\cdot [Z]$ are zero.  Let $\eta=[Y]-[Z]$. Then $h\cdot \eta=0$, so Theorem 
\ref{existence} applies; let $M$ be the module of finite length and finite
projective dimension whose existence is guaranteed by the theorem.  Then if
$\p$ is the graded prime ideal that defines the subscheme $Y$, we have
$$\dime(M)+\dime(A/\p)=0+2=2<3=\dime(A)$$
and
$$\chi(M,A/\p)=\degree(([Y]-[Z])\cdot [Y])=-1.$$
We remark that under an appropriate identification of $\Proj(R)$ with
$\Px^1\times \Px^1$ the ideal $\p$ will be the ideal $(X,Z)$.

Theorem \ref{existence} shows the existence of modules with given intersection
properties but does not show how to construct them.  In joint
work with Greg Piepmeyer, we are developing a method for
explicitly contructing these complexes, and we have
carried it out to the point of writing out a free resolution in
complete detail in the case described in the previous paragraph.

\end{document}